\documentclass[a4paper,11pt]{article}
\usepackage{theorem}
\usepackage{amssymb,amsmath}
\usepackage{latexsym}
\usepackage{color}
\usepackage{txfonts}
\usepackage{multirow}
\usepackage[latin1]{inputenc}

\usepackage{color,graphicx,array}

\theoremheaderfont{\normalfont\bfseries\upshape}
\theorembodyfont{\normalfont\mdseries\slshape}
\newtheorem{theorem}{Theorem}[section]
\newtheorem{proposition}[theorem]{Proposition}
\newtheorem{corollary}[theorem]{Corollary}
\newtheorem{lemma}[theorem]{Lemma}
\theorembodyfont{\normalfont\mdseries\upshape}
\newtheorem{remark}[theorem]{Remark}
\newtheorem{example}[theorem]{Example}
\newtheorem{definition}[theorem]{Definition}

\newcommand{\R}{\mathbb{R}}
\newcommand{\N}{\mathbb{N}}
\newcommand{\I}{\mathcal{K}_{C}}
\newcommand{\Kc}{\mathcal{K}_{C}}
\newcommand{\difgH}{\circleddash_{gH}}

\newcommand{\fu}{\underline{f}}
\newcommand{\fo}{\overline{f}}

\newcommand{\hu}{\underline{h}}
\newcommand{\ho}{\overline{h}}
\newcommand{\continuasintervalares}{\emph{C}\left(\Omega,\I\right)}
\newcommand{\continuasintervalaresab}{\emph{C}\left([a,b],\I\right)}

\begin{document}

\pagestyle{plain} \setcounter{page}{1}
\noindent \textbf{{\LARGE A perturbed collage theorem and its application to inverse interval  integral problems}}

\noindent \textbf{M. Arana--Jim\'enez$^1$, M. I. Berenguer$^2$, D. G\'amez$^2$, A.I Garralda--Guillem$^2$, M. Ruiz Gal\'an$^2$}

\noindent $\mbox{}^1$ Department of Statistics and Operational Research, Faculty of SSCC and Communication, University of C\'adiz, C\'adiz, Spain \\
$\mbox{}^2$ Department of Applied Mathematics, E.T.S. de Ingenier\'ia de Edificaci\'on,  University of Granada, Granada, Spain

\begin{abstract}
	This paper deals with inverse problems subject to imprecise or vague information of some involved data by means of interval-valued functions. To provide interval solutions to the inverse problems we have adopted a perturbed collage-based approach and we have also  introduced a numerical procedure by means of the use of interval bases in a sense. To illustrate the results, and as an application, we have studied the Volterra interval-valued integral equation, and provide some computational examples.
	
\end{abstract}

\noindent {\footnotesize \textbf{2020 Mathematics Subject Classification MSC: }   45Q05, 47S40, 65L10, 65R20.}

\section{Introduction}
\vspace{-0.5cm}

This work deals with some inverse problems by mean of fractal-based methods, and more specifically, collage-based techniques, as well as some related numerical schemes. It
addresses a perturbed version of a well-known and straightforward consequence of the Banach fixed point theorem, the \textsl{collage theorem} \cite{Ervin}, and its application for designing numerical algorithms to solve an inverse problem involving a class of interval integral equations. In order to present the motivation behind that perturbed result, we should first recall that, in its classical form,  the collage theorem asserts that in a complete metric space $(M,d)$, for the unique fixed point $x^\bullet \in M$ of a contractive self-mapping $\Phi : M \longrightarrow M$, the inequality
\[
d(x,x^\bullet) \le \frac{1}{1-c}d(x,\Phi (x)), 
\]
is valid, where $c \in [0,1)$ is the contraction constant and $x$ is any element in $M$. Therefore, a practical problem arises when trying to calculate the so-called \textsl{collage distance}, $d(x,\Phi (x))$: evaluating $\Phi$ at $x$ is not always a feasible calculus from a practical perspective. For instance, this is the case when dealing with integral operators. In our main result, we state the collage theorem in terms of an approximating element $y$ of $\Phi(x)$, which is easy to calculate in a sense, and which satisfies a perturbed inequality along the lines of that which is outlined in the collage theorem.

A question related  to the collage theorem is the collage-based inverse problem: given a target element $x$ in the complete metric space $(M,d)$, a nonempty set $\Lambda$ (typically a nonempty compact subset of a Euclidean real space $\mathbb{R}^N$) and a family $\{\Phi_\lambda\}_{\lambda \in \Lambda}$ of $c_\lambda$-contractive self-mappings ($0\le c_\lambda<1$) on $M$, with respective fixed points $x^\bullet_\lambda$, the aim is to find a parameter $\lambda^* \in \Lambda$ for which the distance $d(x,x^\bullet_{\lambda^*})$ is as small as possible. To this end, and according to the collage theorem, one  must consider (see \cite{Kunze}) the optimization problem
\[
\min_{\lambda \in \Lambda} \frac{1}{1-c_\lambda}d(x,\Phi_\lambda(x)).
\]
Moreover, if $\displaystyle c:=\inf_{\lambda \in \Lambda} c_\lambda  >0$, then we  must instead deal with the nonlinear program 
\[
\min_{\lambda \in \Lambda} d(x,\Phi_\lambda(x)).
\]
For the same practical reasons we have mentioned above, in this work we also focus on with a perturbation of it.

In addition, and as an extension of the previous contents, we have considered another kind of perturbed framework, which is in the presence of interval uncertainty. This is the case when some input (data) is not deterministic, such as, following the example discussed by Rom\'an-Flores and Rojas-Medar \cite{flores-medar02}, in differential equations given for the growing population by means of the Malthus model, in the case of imprecision or vagueness of information from individuals in the initial population. To this matter, and specially on interval arithmetic operations, the reader can find the interval analysis commented on works by Moore \cite{moore66,moore79}. It has a wide range of applications, such as reliable computing, validated numerics, interval problems with differential equations, data envelopment analysis, robotics and so on, which are discussed in several monographs and research papers, in addition to the bibliography therein (see \cite{shwartz16,Stefanini2010, Stefanini2009,stefanini-arana19}).
Alefeld and Mayer \cite{alefeld-mayer10} offer an interesting overview of applications of interval arithmetic, and among these they include the application of the allocation of fixed points under contractive functions, which are point wise valued on intervals. Furthermore, and in relation to our objectives, in  recent literature we  find definitions of integral for interval-valued functions, such as those given by Aubin and Cellina \cite{Aubin-Cellina84}, Wu and Gong \cite{WuGong00}, and Stefanini and Bede \cite{Stefanini2009}, among others, which will be referred to and discussed in this text. The integral for interval-valued functions provided by Stefanini and Bede \cite{Stefanini2009} will be a useful tool for us to use as a model and  tooperate inverse problems subject to interval uncertainty, such as the Volterra interval integral equation, as proposed in the present manuscript.

 This leads us to consider interval problems, and we have focused on those of  an integral-type. In this context, the set $C(\Omega, \mathcal{K}_c)$ of continuous functions from a compact space $\Omega$ into the metric space $\mathcal{K}_c$ of all compact real intervals (the Hausdorff metric), becomes a complete metric space when endowed with an adequate uniform metric. Therefore, the interval integral equation under study permits the described treatment of the corresponding inverse problem. Another key factor is the way in which we construct the perturbations in this metric space. We approximate a continuous interval-valued function by means of a Schauder basis in an associated space of continuous functions, which generates, adequately way, easy-to-calculate continuous interval-valued functions to evaluate the involved operators $\Phi_\lambda$. So, for a given function $X \in C(\Omega, \mathcal{K}_c)$, we find another one $Y \in C(\Omega, \mathcal{K}_c)$, as close to $\Phi_\lambda(X)$ as we want, in such a way that we can explicitly determine $Y$.

The structure of this paper is as follows. In Section \ref{sec:perturbedcollage} we derive the perturbed collage theorem as a consequence of the (Caccioppoli version of the) Banach fixed point theorem and in Section \ref{sec:preliminaries} we recall some facts about interval calculus. In order to performthe calculations involved in the inverse problem associated with a certain interval integral equation  effectively,  in Section \ref{sec:approximation} we have designed an easy-to-calculate procedure derived from the use of some Schauder bases in certain spaces of continuous functions. Finally, Section \ref{secinverse} deals with the study of the above-mentioned inverse problem, that of the Volterra interval-valued integral equation

\begin{equation} \label{eq:volterra}
X(t)=G(t)+\int_{a}^{t} K(t,s,X(s))ds, \qquad (t \in [a,b]),
\end{equation}
where $G\in \emph{C}([a,b],\mathcal{K}_c)$, $K \in \emph{C}([a,b] \times [a,b]\times \mathcal{K}_c,\mathcal{K}_c)$ are assumed to be known interval functions, and $X \in \emph{C}([a,b],\I)$ is the unknown interval-valued function to be determined. Finally, in Section \ref{sec:examples} we provide some examples.

\section{A perturbed collage theorem}\label{sec:perturbedcollage}
\vspace{-0.5cm}

In this paper we have developed a collage-based result, but as previosly mentioned in the Introduction, in practical situations we can not explicitly determine the collage distance, only an approximation of it. In order to state the corner stone of these ideas, we begin by establishing a collage theorem of perturbed character, which follows from the generalization of R. Caccioppoli of the Banach fixed point theorem (see, for instance, \cite[Theorem 2.3]{alm-ans-kha14}). For a self-mapping $\Phi$ on a nonempty set $M$ and $x \in M$, $\{\Phi^n(x)\}_{n \ge 1}$ denotes the sequence of iterates of $\Phi$ generated by $x$. Then, that extension of Banach's fixed point theorem reads as follows: assume that $(M,d)$ is a complete metric space, $\Phi : M \longrightarrow M$ such that      
\[
x_1,x_2 \in M \ \Rightarrow \ d(\Phi^n(x_1),\Phi^n(x_2)) \le \alpha_n d(x_1,x_2)
\]
 and $\{\alpha_n\}_{n \ge 1} \subset \mathbb{R}_+$ is a sequence such that the series $\displaystyle \sum_{n \ge 1} \alpha_n $ is convergent. Then, $\Phi$ has a unique fixed point $x^\bullet \in M$ and if $x \in M$, then  for all $n\ge 1$ there holds that
\[
d(\Phi^n(x),x^\bullet) \le \sum_{k=n}^\infty \alpha_k d(\Phi(x),x),
\]
and, in particular, $\displaystyle \lim_{n \to \infty} \Phi^n(x)=x^\bullet$.

\begin{theorem}\label{th:perturbedcollage}
Let $(M,d)$ be a complete metric space, $\Phi : M \longrightarrow M$ and $\{\alpha_n\}_{n \ge 1}$ be a sequence of nonnegative real numbers such that $\displaystyle \sum_{n=1}^\infty \alpha_n < + \infty$ and
\[
x_1,x_2 \in M \ \Rightarrow \ d(\Phi^n(x_1),\Phi^n(x_2)) \le \alpha_n d(x_1,x_2).
\]
If in addition $x^\bullet$ is the unique fixed point of $\Phi$, and $\varepsilon >0$ and $x, y \in M$ satisfy $d(\Phi(x),y)<\varepsilon$, then for each $n\ge 1$ with $\alpha_n <1$ the inequality
\begin{equation}\label{eq:desigualdad}
d(x,x^\bullet) \le \frac{\displaystyle \sum_{k=1}^{n-1}\alpha_k}{1-\alpha_n}(d(x,y)+\varepsilon)
\end{equation}
is valid.
\end{theorem}

\noindent \textit{Proof.} First of all, we shoul note that the above-mentioned extension of the Banach fixed point  theorem, \cite[Theorem 2.3]{alm-ans-kha14}, guarantees the existence of a unique fixed point $x^\bullet$ of $\Phi$. So, to conclude this proof, we should fix $\varepsilon >0$ and $x,y \in M$ satisfying $d(\Phi(x),y)<\varepsilon$, and consider  $n \ge 1$ so that $\alpha_n < 1$ (the existence of such an $n \ge 1$, indeed, that of all them except perhaps a finite number, follows from the convergence of the series $\displaystyle \sum_{n \ge 1} \alpha_n$). Therefore,
\[
\begin{array}{rl}
d(x,x^\bullet) & \le d(x,\Phi^n(x))+d(\Phi^n(x),x^\bullet)   \\
               & = d(x,\Phi^n(x))+d(\Phi^n(x),\Phi^n(x^\bullet))   \\
               & \le d(x,\Phi^n(x))+\alpha_n d(x,x^\bullet),   
\end{array}
\]
and thus
\[
d(x,x^\bullet) \le \frac{1}{1-\alpha_n} d(x,\Phi^n(x)).
\]
As a consequence, we arrive at the announced inequality, since
\[
\begin{array}{rl}
d(x,x^\bullet) & \le \displaystyle \frac{1}{1-\alpha_n} d(x,\Phi^n(x))   \\
               & \le \displaystyle \frac{1}{1-\alpha_n} \sum_{k=0}^{n-1} d(\Phi^k(x)),\Phi^k(\Phi (x))   \\
               & \le \displaystyle \frac{\displaystyle\sum_{k=0}^{n-1} \alpha_k}{1-\alpha_n}d(\Phi(x),x)  \\
               & \le \displaystyle \frac{\displaystyle\sum_{k=0}^{n-1}\alpha_k}{1-\alpha_n}(d(\Phi(x),y)+d(y,x))   \\  
               & \le \displaystyle \frac{\displaystyle\sum_{k=0}^{n-1}\alpha_k}{1-\alpha_n}(\varepsilon+d(y,x)). 
\end{array}
\]
\hfill $\square$

\section{Some elementary notions and facts about interval calculus}\label{sec:preliminaries}
\vspace{-0.5cm}

We denote by $\mathcal{K}_{C}$ the family of all bounded closed intervals in $\mathbb{R}$, i.e.,
\begin{equation*}
\mathcal{K}_{C}=\left\{ \left[ \underline{a},\overline{a}\right] : \ \underline{a},\overline{a}\in \mathbb{R}\mbox{ and
}\underline{a}\leq \overline{a}\right\}.
\end{equation*}%
From herein, for the sake of simplicity, we refer to any element in $\Kc$ as an interval.

Given two intervals
$A=[\underline{a},\overline{a}],\;B=[\underline{b},\overline{b}]$, we have the following classical operations of the sum
\[
A+B=[\underline{a}+\underline{b},\overline{a}+\overline{b}]
\]
and the multiplication by scalars 
\[
\tau A=\{\tau a :\  a\in A\}=\left\{
	\begin{array}{ll}
	[\tau \underline{a},\tau \overline{a}], & \mbox{ if } \tau\ge 0,\\
	
	[\tau\overline{a},\tau\underline{a}], & \mbox{ if } \tau\leq 0.
	\end{array}\right.
\]
With respect to the difference of two intervals, several definitions exist in literature. One of the most popular is the generalized Hukuhara difference ($gH$-difference, for short). The $gH$-difference of two intervals $A$ and $B$, which we recall from \cite{markov79,Stefanini2010,Stefanini2009}, is as follows:
\begin{equation*}
A\ominus_{gH}B=C \Longleftrightarrow
\left\{
\begin{array}{ll}
&(a)\ A= B+C,\\
or&(b)\ B = A + (-1)C.
\end{array}
\right.
\end{equation*}
Note that the difference of an interval and itself is zero, that is,  $A\ominus_{gH}A=[0,0]$.
Furthermore, the $gH$-difference of two intervals always exists and is equal to
\begin{equation*}
A\ominus_{gH}B = [\min\{\underline{a}-\underline{b}, \overline{a}-\overline{b}\}, \max\{\underline{a}-\underline{b}, \overline{a}-\overline{b}\}].
\end{equation*}
We refer to Markov \cite{markov79}, Moore \cite{moore66,moore79} and Alefeld and Herzberger \cite{alefeld83} for further details on the topic of interval
analysis. 

We also recall the \textsl{Pompeiu-Hausdorff} distance $D$ on $\mathcal{K}_{C}$  which  is defined by
\begin{equation*}
D(A,B) = \max\left[\max_{a\in A}d(a,B),\max_{b\in B}d(b,A)\right]
\end{equation*}
with $d(a,B)=\min_{b\in B} |a-b|$. It is well-known (see \cite{Stefanini2009}) that $D(A,B)=\left\Vert A\ominus
_{gH}B\right\Vert = \max\{|\underline{a}-\underline{b}|, |\overline{a}-\overline{b}|\}$ where, for $C\in \mathcal{K}_{C}$, $\left\Vert
C\right\Vert =\max \{\left\vert c\right\vert :\ c\in C\}$, and that $(\Kc,D)$ is a complete and separable metric space. 

Finally, in this paper we consider the definition of the integral for an interval-valued function as used by L. Stefanini and B. Bede \cite{Stefanini2009}, as follows.
\begin{definition}\label{defi: integral}
	Given an interval-valued function $f:[a,b]\rightarrow \I$, with $f(t)=[\underline{f}(t), \overline{f}(t)]$, $f$ is said to be integrable if and only if the endpoint functions $\underline{f}$ and $\overline{f}$ are integrable, and the integral of $f$ over $[a,b]$ is defined as
	\begin{equation}\label{def: integral}
	\int_{a}^{b} f(t)dt = \left[ \int_{a}^{b}\underline{f}(t)dt, \int_{a}^{b} \overline{f}(t)dt \right]\in\Kc.
	\end{equation}
\end{definition}

\begin{remark}
	The previous definition is a natural extension of the classic integral, and fits the aim of this paper. However, we would like to point out that other definitions of integral for inteval-valued functions exist in literature. For instance, Aubin and Cellina \cite{Aubin-Cellina84} have defined the Aunmann integral of $f$ over $[a,b]$ as $\int_{a}^{b} f(t)dt=\left\{ \int_{a}^{b} h(t)dt:h\in S(f)  \right\}$, where $S(f)$ is the set of all the integrable selectors of $f$, that is, $S(f)=\left\{ h:[a,b]\rightarrow\R: h \mbox{ integrable, } h(t)\in f(t),\forall t\in[a,b]\right\} $. Bede and Gal \cite{Bede-Gal05} have proved that if $f$ is a measurable and integrable bounded interval-valued function, then the endpoint functions $\underline{f}$ and $\overline{f}$ are integrable, and (\ref{def: integral}) is fulfilled. On the other hand, Wu and Gong \cite{WuGong00} introduced the Henstock integral for interval-valued functions as an extension of the Henstock integral in real-valued functions. In fact, the endpoint functions $\underline{f}$ and $\overline{f}$ are Henstock integrable if and only if the interval-valued function $f$ is Henstock integrable, and its integral coincides with the interval defined by the Henstock integral of the endpoint functions (\cite{WuGong00}, Theorem 2.1). In this way, under certain hypotheses, the Henstock integral of $f$  verifies (\ref{def: integral}).
\end{remark}

\section{Approximation of a continuous interval-valued function}\label{sec:approximation} 
\vspace{-0.5cm}

In this section we focus on obtaining an approximation of a continuous interval-valued function in terms of a sequence of functions, simple in a sense, which will be essential in the development of our numerical method to solve interval integral inverse problems by means of the perturbed collage theorem, Theorem \ref{th:perturbedcollage}.
 
Given $\Omega$ a compact topological space, we will denote
 by $\emph{C}\left(\Omega\right)$ the Banach space of all continuous  real valued functions defined on $\Omega$ with its usual $\max$ norm, $\Vert \cdot \Vert_\infty$, and  by $\continuasintervalares$ the set of all continuous functions from $\Omega$ into $\Kc$ endowed with the distance 
\begin{equation}\label{eq:distance}
H(f,h):=\sup_{\omega\in\Omega} \max\{|\underline{h}(\omega)-\fu(\omega)|, |\overline{h}(\omega)-\fo(\omega)|\}
\end{equation}
with  $f(\omega)=[\fu(\omega),\fo(\omega)]$ and $h(\omega)=[\underline{h}(\omega),\overline{h}(\omega)]$ in $\continuasintervalares$.

It is very easy to check that $f=[\fu,\fo]:\Omega \rightarrow \Kc$ is continuous if and only if $\fu$  and $\fo$ are continuous. As a direct consequence, if $T_\Omega$ is the nonempty set
$$
T_\Omega:=\{(\fu,\fo) \in \emph{C}\left(\Omega\right)\times  \emph{C}\left(\Omega\right): \fu\le \fo\}
$$
with the metric
$$ d((\fu,\fo), (\hu,\ho)):= \max\{\Vert\fu-\hu\Vert_{\infty}, \Vert \fo-\ho\Vert_{\infty} \},$$
then, the mapping
$$
S: (\continuasintervalares,H) \rightarrow (T_\Omega,d)
$$
defined by
$$S([\fu,\fo]):=(\fu,\fo)$$
is an isometric bijection. In particular, $(\continuasintervalares,H)$ is a complete and separable metric space.

In order to obtain the aforementioned approximation for a continuous interval-valued function, we make use of an adequate Schauder basis in a certain Banach space, motivated by the succesful use in the scalar-valued case (see \cite{Berenguer2017},  \cite{Berenguer2011}, \cite{Berenguer2010} and \cite{Berenguer2009}). We should recall that, in general, given a Banach space $E$, a sequence $\{f_n\}_{n\ge 1}$  of elements of $E$ is said to be a \textsl{Schauder basis} of $E$ if, for every $z\in E$, there is a unique sequence
$\{\alpha_n\}_{n\ge 1}$ of scalars such that $z=\sum_{n\ge 1} \alpha_n f_n$. A Schauder basis gives rise to the canonical sequence of (continuous and linear) associated \textsl{projections} $\Pi_m:E\rightarrow E$, $\Pi_m (\sum_{n\ge 1}
\alpha_n f_n) := \sum_{k=1}^m \alpha_k f_k$ (see \cite{gelbaum} and \cite{semadeni1982}).

In the following result we introduce an approximation of a continuous interval-valued function in  the metric space $\emph{C}(\Omega, \Kc)$  in terms of a Schauder basis in the Banach space $\emph{C}(\Omega)$, which will be essential for obtaining the algorithm to solve the inverse problem for the Volterra interval integral equation.

\begin{proposition}\label{pr:schauder}  Let $\Omega$ be a topological compact space,  $\{f_n\}_{n\ge 1}$ be a Schauder basis of $\emph{C}(\Omega)$ and  let  $\{ \Pi_n\}_{n \ge 1}$  be the associated  sequence of projections such that 
\begin{enumerate}
\item[a)] if $\omega \in \Omega$ and $n \ge 1$, then $f_n(\omega)\geq 0$, and
\item[b)] if $g \in C(\Omega)$, $g \ge 0$, and $n \ge 1$, then $\Pi_n(g)\ge 0$.
\end{enumerate}
Thus, given $f \in \emph{C}(\Omega, \Kc)$ and $\varepsilon>0$ there exists $n\ge 1$ such that 
$$H\left( f, P_n(f)\right)< \varepsilon,$$
where 
$$
P_n(f)(\omega)=\sum_{k=1}^n \alpha_k \phi_k (\omega) + \sum_{\substack{k=1 \\ \beta_k-\alpha_k \ge 0 }}^n (\beta_k-\alpha_k) \psi_k(\omega)\difgH \left(\sum_{\substack{k=1 \\ \beta_k-\alpha_k < 0 }}^n  \vert\beta_k-\alpha_k\vert  \psi_k(\omega) \right)
$$ 
with $\phi_k(\omega)= \left[f_k(\omega),f_k(\omega)\right]$, $\psi_k(\omega)=\left[0, f_k(\omega)   \right]$ and certain real numbers $\alpha_1,\alpha_2,\dots, \alpha_n$ and $\beta_1, \beta_2,\dots,\beta_n$.
\end{proposition}

\noindent \textit{Proof}. If $f=[\fu,\fo]\in \continuasintervalares$, the fact  that $\{f_n\}_{n \ge 1}$ is a Schauder basis of $C(\Omega)$ implies the existence of sequences of scalars $\{\alpha_n\}_{n \ge 1}$ and $\{\beta_n\}_{n \ge 1}$ such that $\fu = \sum_{n\ge 1} \alpha_n f_n$ and $\fo= \sum_{n\ge 1} \beta_n f_n$. Therefore, we define 
$$
P_n(f) (\omega):=\left[  \sum_{k=1}^n \alpha_k f_k(\omega),\sum_{k=1}^n \beta_k f_k(\omega)   \right], \qquad (n \ge 1, \ \omega \in \Omega),
$$
and, according to the definition \eqref{eq:distance} of the distance $H$, given $\varepsilon >0$ we can guarantee the existence of $n\ge 1$ such that $H(f, P_n(f))<\varepsilon$. Finally, taking into account that if $a,b,\alpha,\beta\in\R^{+}_0$ with $a\alpha-b\beta>0$ , then
\[ 
[0,a\alpha-b\beta]=a[0,\alpha]\difgH b[0,\beta],
\]
it follows that
$$
\begin{array}{rl}
P_n(f) (\omega) & = \displaystyle \left[\sum_{k=1}^n \alpha_k f_k(\omega), \sum_{k=1}^n \alpha_k f_k(\omega)\right]+ \left[0,\sum_{k=1}^n (\beta_k-\alpha_k) f_k(\omega)\right]   \\
     & = \displaystyle \sum_{k=1}^n \alpha_k  \left[f_k(\omega),f_k(\omega)\right]+\left[0,\sum_{\substack{k=1 \\ \beta_k-\alpha_k \ge 0 }}^n (\beta_k-\alpha_k) f_k(\omega)   \right]\difgH \left[0,\sum_{\substack{k=1 \\ \beta_k-\alpha_k < 0 }}^n \vert\beta_k-\alpha_k\vert f_k(\omega)   \right] \\   
     & = \displaystyle \sum_{k=1}^n \alpha_k  \left[f_k(\omega),f_k(\omega)\right]+ \sum_{\substack{k=1 \\ \beta_k-\alpha_k \ge  0 }}^n (\beta_k-\alpha_k) \left[0, f_k(\omega)   \right]\difgH \left(\sum_{\substack{k=1 \\ \beta_k-\alpha_k < 0 }}^n  \vert\beta_k-\alpha_k\vert  \left[0, f_k(\omega)   \right]\right)  \\
     & \displaystyle = \sum_{k=1}^n \alpha_k \phi_k (\omega) + \sum_{\substack{k=1 \\ \beta_k-\alpha_k \ge  0 }}^n (\beta_k-\alpha_k) \psi_k(w)\difgH \left(\sum_{\substack{k=1 \\ \beta_k-\alpha_k < 0 }}^n  \vert\beta_k-\alpha_k\vert  \psi_k(\omega) \right), 
\end{array}
$$ 
where  $\phi_k(\omega)= \left[f_k(\omega),f_k(\omega)\right]$, $\psi_k(\omega)=\left[0, f_k(\omega)   \right]$.
\hfill $\square$

The following result provides us with an procedure to determine  the integral of $P_n(f)$ in Proposition \ref{pr:schauder}:

\begin{proposition}\label{prop: properties of integral}
	(i) Given $h_k:[a,b]\rightarrow \Kc$ integrable interval-valued functions, for $k=1,\dots, n$, then $\sum_{k=1}^{n} h_k$ is integrable, and $\int_{a}^{b}\sum_{k=1}^{n} h_k(t)dt=\sum_{k=1}^{n}\int_{a}^{b} h_k(t)dt$.\\
	(ii) Given $f,\; h:[a,b]\rightarrow \Kc$ integrable interval-valued functions, and define $(h\ominus_{gH} f)(t)=h(t)\ominus_{gH} f(t)$, then $h\ominus_{gH} f$ is integrable. Furthermore, if $f(t)=[0,\overline{f}(t)]$ and $h(t)=[0,\overline{h}(t)]$, with $\overline{h}(t)\ge \overline{f}(t)$, then  $\int_{a}^{b}(h(t)\ominus_{gH} f(t))dt=\int_{a}^{b}h(t)dt\ominus_{gH} \int_{a}^{b}f(t)dt$.
\end{proposition}

\noindent \textit{Proof}.
		(i) The proof is straightforward, since the summation of $h_k$, $k=1,\dots,n$, is equivalent to the summantion of their lower and upper endpoints, which are integrable by Definition \ref{defi: integral}.\\
		(ii) By definition of the $gH$-difference, we have that 
		$$h(t)\ominus_{gH} f(t)=\left[ \mbox{min}\{\underline{h}(t)-\underline{f}(t),\overline{h}(t)-\overline{f}(t)\},\mbox{max}\{\underline{h}(t)-\underline{f}(t),\overline{h}(t)-\overline{f}(t)\}  \right].$$
		Since $f$ and $h$ are integrable, then $\underline{h}, \;\underline{f},\; \overline{h},\;\overline{f}$ are integrable. The latest implies that $\underline{h}-\underline{f}$ and $\overline{h}-\overline{f}$ are integrable. Then, $\mbox{min}\{\underline{h}-\underline{f},\overline{h}-\overline{f}\}$ and $\mbox{max}\{\underline{h}-\underline{f},\overline{h}-\overline{f}\}$ are integrable, therefore $h\ominus_{gH} f$ is integrable. Furthermore if $f(t)=[0,\overline{f}(t)]$ and $h(t)=[0,\overline{h}(t)]$ with $\overline{h}(t)\ge \overline{f}(t)$, then $\int_{a}^{b}(h(t)\ominus_{gH} f(t))dt=\left[0,\int_{a}^{b}(\overline{h}(t)-\overline{f}(t))dt  \right]=
		\int_{a}^{b}h(t)dt\ominus_{gH} \int_{a}^{b}f(t)dt$.\hfill $\square$

Now we are in a position to calculate the integral of $P_n(f)$.

\begin{proposition}\label{prop: computational operations}
	Consider the interval-valued functions and the hypotheses and equations given in Proposition \ref{pr:schauder} with $\Omega:=[a,b]$. Then,
\[
	\int_{a}^{b}P_n(f)(t)dt=\sum_{k=1}^n \alpha_k \int_{a}^{b}\phi_k (t)dt + \sum_{\substack{k=1 \\ \beta_k-\alpha_k \ge 0 }}^n (\beta_k-\alpha_k)\int_{a}^{b} \psi_k(t)dt\difgH \left(\sum_{\substack{k=1 \\ \beta_k-\alpha_k < 0 }}^n  \vert\beta_k-\alpha_k\vert  \int_{a}^{b}\psi_k(t)dt \right).
\]
\end{proposition}

\noindent \textit{Proof}.  The interval-valued functions $\phi_k$  and $\psi_k$ verify the hypothesis required in (i) and (ii) in Proposition \ref{prop: properties of integral}. Therefore, in view of Proposition \ref{pr:schauder}, the validity of the announced equality follows. \hfill $\square$

The following technical result together with Proposition \ref{prop: properties of integral} and Propostion  \ref{prop: computational operations}, allow us to check that $\int_{a}^{b}P_n(f)(t)dt$ is close to  $\int_{a}^{b}f(t)dt$ in the sense of the distance $D$.

\begin{lemma}\label{le:modular}
If  $f, h \in \emph{C}([a,b],\mathcal{K}_c)$, then
\[
D \left( \int_a^b f(s)ds,\int_a^b h(s)ds \right) \le (b-a) H(f,h). 
\]
\end{lemma}

\noindent \textit{Proof}. If $f(s)=[\underline{f}(s),\overline{f}(s)]$ and  $h(s)=[\underline{h}(s),\overline{h}(s)]$, it suffices to follow this chain of inequalities:
\[
\begin{array}{rl}
\displaystyle D \left( \int_a^b f(s)ds,\int_a^b h(s)ds \right) & = \displaystyle \max \left\{ \left| \int_a^b \underline{f}(s)ds-\int_a^b\underline{h}(s)ds \right| ,  \left| \int_a^b \overline{f}(s)ds-\int_a^b\overline{h}(s)ds \right|  \right\}  \\
     & \displaystyle \le \max \left\{ \int_a^b | \underline{f}(s)-\underline{h}(s)|ds  ,  \int_a^b | \overline{f}(s)-\overline{h}(s)| ds   \right\}   \\
          & \displaystyle \le \int_a^b \max \left\{ | \underline{f}(s)-\underline{h}(s)|  , | \overline{f}(s)-\overline{h}(s)|    \right\} ds   \\
          & = \displaystyle \int_a^b D(f(s),h(s))ds \\
& \leq  \displaystyle \int_a^b H(f,h) ds \\
& =(b-a) H(f,h).
\end{array}
\]
\hfill$\square$

\section{Inverse problem}\label{secinverse}
\vspace{-0.5cm}

Since, in this work, we adopt a collage-based approach to deal with an inverse problem related to the Volterra interval integral equation \eqref{eq:volterra}, we must first consider its fixed-point-treatment, that is, the well-known fact that, under suitable assumptions of continuity and lipschitzianity, the equation \eqref{eq:volterra} admits a unique solution, and that it is possible to give this in terms of the limit of the sequence of iterates of an integral operator at a continuous interval-valued function, since such an interval integral equation is obviously equivalent to determine a fixed point of an adequate integral operator. This follows from the generalized Banach fixed point theorem (\cite[Theorem 2.3]{alm-ans-kha14}).

\begin{proposition}\label{pr:existence}
Assuming $G \in \emph{C}([a,b],\mathcal{K}_c)$ and $K \in \emph{C}([a,b]^2\times \mathcal{K}_c,\mathcal{K}_c)$, in such a way that, for some $L>0$ there holds that
\[
\left.
\begin{array}{c}
A,B \in \mathcal{K}_c \\
a \le t,s \le b
\end{array}
\right\}
\ \Rightarrow \
D(K(t,s,A),K(t,s,B)) \le L\  D(A,B). 
\] 
Then, the self-operator $\Phi : \emph{C}([a,b],\mathcal{K}_c) \longrightarrow \emph{C}([a,b],\mathcal{K}_c)$ defined at each $X \in \emph{C}([a,b],\mathcal{K}_c)$ as
\begin{equation}\label{equVol2}
\Phi(X)(t):=G(t)+\int_{a}^{t} K(t,s,X(s))ds, \qquad (t \in [a,b])
\end{equation}
admits a unique fixed point $X^\bullet$. Furthermore, if $X \in \emph{C}([a,b],\mathcal{K}_c)$ and $n \ge 1$, then
\[
H(\Phi^n(X),X^\bullet) \le \left( \sum_{k=n}^\infty \frac{L^k(b-a)^k}{k!}\right) H (\Phi(X),X),
\]
and so,
\[
\lim_{n \to \infty} \Phi^n(X)=X^\bullet.
\]
\end{proposition}

In particular, and according to the perturbed collage theorem, Theorem \ref{th:perturbedcollage}, we arrive at:

\begin{corollary}\label{co:perturbedcollage}
Let  $G \in \emph{C}([a,b],\mathcal{K}_c)$, $K \in \emph{C}([a,b]^2\times \mathcal{K}_c,\mathcal{K}_c)$ and $L>0$ such that
\[
\left.
\begin{array}{c}
A,B \in \mathcal{K}_c \\
a \le t,s \le b
\end{array}
\right\}
\ \Rightarrow \
D(K(t,s,A),K(t,s,B)) \le L\  D(A,B). 
\] 
Let $X^\bullet \in \emph{C}([a,b],\mathcal{K}_c)$ be the unique solution of the Volterra interval integral equation \eqref{eq:volterra} and suppose that $\varepsilon >0$ and $X,Y \in \emph{C}([a,b],\mathcal{K}_c)$ satisfy 
\[
H \left( G(\cdot)+\int_{a}^{(\cdot)} K(\cdot,s,X(s))ds,Y(\cdot) \right) < \varepsilon .
\]
Then
\[
H(X,X^\bullet) \le e^{L(b-a)}(H(X,Y)+\varepsilon).
\]
\end{corollary}

With the idea in mind of stating the previosly inverse problem in a precise way, now we can focus on the following fact: In the perturbed collage theorem, Theorem \ref{th:perturbedcollage}, we can replace the image by $\Phi$ of an element $x \in M$ by another $y \in M$ where, for a given $\varepsilon >0$, $d(x,y)<\varepsilon$. When dealing with a certain inverse problem involving the Volterra interval integral equation \eqref{eq:volterra}, our aim is to find, for a given $X \in C([a,b],\mathcal{K}_c)$ and an $\varepsilon >0$, a continuous function $Y \in C([a,b],\mathcal{K}_c)$ such that both $H(\Phi(X),Y)<\varepsilon$ and that $Y$ is easy-to-calculate, unlike $\Phi(X)$ in general. We do this in  Lemma \ref{le:preparacion} in a constructive way, which will be essential for addressing the inverse problem.

\begin{lemma}\label{le:preparacion}
Let $\{f_n\}_{n \ge 1}$ be a Schauder basis in $\emph{C}([a,b]^2)$ with sequence of associated projections $\{\Pi_n\}_{n \ge 1}$ and satisfying the hypotheses in Proposition \ref{pr:schauder}, $X \in \emph{C}([a,b],\mathcal{K}_c)$, $\Phi$ be defined as in \eqref{equVol2} and $\varepsilon >0$. Therefore, there exists $n \ge 1$ in such a way that for $P_n$ defined as in Proposition \ref{pr:schauder} and $Y:[a,b] \longrightarrow\mathcal{K}_c$ is the continuous interval-valued function
\[
Y(\cdot):=G(\cdot)+\int_a^{(\cdot)}P_n(K(\cdot,s,X(s))ds,
\]
we have that
\[
H(\Phi(X),Y)<\varepsilon.
\]
\end{lemma} 

\noindent \textit{Proof}. If $Z \in \emph{C}([a,b]^2,\mathcal{K}_c)$ is defined at each $a \le t,s \le b$ by
\begin{equation}\label{eq:z}
Z(t,s):=K(t,s,X(s)),
\end{equation}
according to Proposition \ref{pr:schauder}, let $n$ be a positive integer with
\[
H(Z,P_n(Z))<\frac{\varepsilon}{b-a}.
\] 
Then, taking into account the invariance of $D$ by translations,  Lemma \ref{le:modular} and \eqref{eq:z}, we arrive at
\[
\begin{array}{rl}
H(\Phi(X),Y) & = \displaystyle \sup_{t \in [a,b]} D\left( G(t)+\int_a^t K(t,s,X(s))ds, G(t)+\int_a^t P_n(K(t,s,X(s)))ds \right)  \\
     & \displaystyle =  \sup_{t \in [a,b]} D\left( \int_a^t K(t,s,X(s))ds, \int_a^t P_n(K(t,s,X(s)))ds \right)  \\
     & \le \displaystyle \sup_{t \in [a,b]} (t-a)  H (K(t,\cdot,X(\cdot)), P_n(K(t,\cdot,X(\cdot))))    \\
     & \le  (b-a)    H(Z,P_n(Z)) \\
& <  \varepsilon.
\end{array}
\]
\hfill$\square$

Finally, we present the general scheme for solving an inverse problem via the perturbed collage theorem, Theorem \ref{th:perturbedcollage}. The idea is to consider a complete metric space $(M,d)$, a nonempty index set $\Lambda$ --usually a compact subset of $\mathbb{R}^N$--, a family of contractive self-mappings $\Phi_\lambda: M \longrightarrow M$ or, at least,those which satisfy the conditions of the perturbed collage theorem, Theorem \ref{th:perturbedcollage}, with a unique fixed point $x^\bullet_\lambda$, a target element $x \in M$ and an $\varepsilon >0$. Then, for any $\lambda \in \Lambda$ we calculate an element $y_\lambda \in M,$ in a suitable way, such that $d(\Phi_\lambda(x),y_\lambda)$ is easily computable and
\[
d(\Phi_\lambda(x),y_\lambda)< \varepsilon,
\]
and then we determine that
\begin{equation}\label{eq:derecha}
\lambda^* \in \underset{\lambda \in \Lambda}{\mathrm{arg min}} \ d(y_\lambda,x),
\end{equation}
for which, in view of \eqref{eq:desigualdad}, it holds that
\begin{equation}\label{eq:intermedio}
d(x,x_\lambda^\bullet) \le \frac{\displaystyle \sum_{k=1}^{n-1}\alpha_k^{(\lambda)}}{1-\alpha_{n_\lambda}^{(\lambda)}}(d(x,y_\lambda)+\varepsilon),
\end{equation}
where, for each $\lambda \in \Lambda$, $\{\alpha_{n_\lambda}\}_{n \ge 1}$ is a sequence in $\mathbb{R}_+$ with $\displaystyle \sum_{n=1}^\infty \alpha_{n_\lambda} < + \infty$ and $n_\lambda\ge 1$ satisfies $\alpha_{n_\lambda} <1$. If, in addition, we assume the stability condition
\begin{equation}\label{eq:stability}
\rho:=\sup \left\{ \frac{\displaystyle \sum_{k=1}^{n-1}\alpha_k^{(\lambda)}}{1-\alpha_{n_\lambda}^{(\lambda)}} : \ \lambda \in \Lambda \right\} < \infty,
\end{equation}
then  \eqref{eq:intermedio} yields
\[
d(x,x_\lambda^\bullet) \le \rho (d(x,y_\lambda)+\varepsilon),
\]
and so we can consider problem \eqref{eq:derecha} instead of
\[
\lambda^* \in \underset{\lambda \in \Lambda}{\mathrm{arg min}} \ d(x,x_\lambda^\bullet),
\] 
since in general we can not determine each fixed point $x_\lambda^\bullet$ or, if it is possible to obtain an approximation of it, such a calculation has a very high computational cost. The condition \eqref{eq:stability} is quite familiar in some particular cases. For instance, when the self-mappings $\Phi_\lambda$ are $c_\lambda$-contractive, with $0 \le c_\lambda < 1$, it is equivalent to the well-known assumption $\displaystyle \sup_{\lambda \in \Lambda} c_\lambda <1$ in a typical collage-based approach to inverse problems (see for instance \cite{Berenguer2016}, \cite{Kunze} and \cite{KunzeTorre}).

Now we can address the inverse problem related to the Volterra interval integral equation \eqref{eq:volterra}. So we consider $a,b \in \mathbb{R}$ with $a<b$, a nonempty set $\Lambda$ and, for each $\lambda \in \Lambda$, two interval-valued functions $G_\lambda \in \emph{C}([a,b],\mathcal{K}_c)$ and $K_\lambda \in \emph{C}([a,b]^2\times \mathcal{K}_c,\mathcal{K}_c)$, and an $L_\lambda>0$ such that
\[
\left.
\begin{array}{c}
A,B \in \mathcal{K}_c \\
a \le t,s \le b
\end{array}
\right\}
\ \Rightarrow \
D(K_\lambda(t,s,A),K_\lambda(t,s,B)) \le L_\lambda \  D(A,B). 
\] 
We should also consider the self-operator $\Phi_\lambda : \emph{C}([a,b],\mathcal{K}_c) \longrightarrow C([a,b],\mathcal{K}_c)$ defined for each $X \in \emph{C}([a,b],\mathcal{K}_c)$ by
\begin{equation}\label{equVol2b}
\Phi_\lambda(X)(t):=G_\lambda(t)+\int_{a}^{t} K_\lambda(t,s,X(s))ds, \qquad (t \in [a,b]),
\end{equation}
and its unique fixed point $X_\lambda^\bullet$. 
Then, in view of  Proposition \ref{pr:existence}, the stability condition \eqref{eq:stability} is valid as soon as 
\[
\sup_{\lambda \in \Lambda} L_\lambda < \infty
\]
and therefore, under this hypothesis, the preceding reasoning applies.

\section{Numerical examples}\label{sec:examples}
\vspace{-0.5cm}

In order to illustrate the behaviour of the inverse problem, we have developed some examples. Given $\Lambda$ a compact subset of $\mathbb{R}^N$, for $\lambda \in \Lambda$, given 
$G_\lambda\in\continuasintervalaresab$, and $K_\lambda:[a,b]\times[a,b]\times \I \rightarrow \I$ we consider the problem 
of finding $X\in \continuasintervalaresab$ such that
\[
X(t)=G_{\lambda}(t)+\int_{a}^{t} K_\lambda(t,s,X(s))ds, \qquad (t \in [a,b]).
\]

The inverse problem starts with a target element $\tilde{X}$. This element is obtained as follows: first of all, we fix values $\lambda_0\in \Lambda$ and we obtain a numerical approximation of the solution of the equation
\begin{equation}\label{eq:vol}
X(t)=G_{\lambda_0}(t)+\int_{a}^{t} K_{\lambda_0}(t,s,X(s))ds, \qquad (t \in [a,b]).
\end{equation}
This numerical approximation is an approximation of the fixed point of the self-operator given in \eqref{equVol2b}. To obtain $\tilde{X}$ we start with an initial $X_0\in \continuasintervalaresab$ and for $j=1,2,...$ we calculate  for the  chosen $n\in \mathbb{N}$ 
$$
X_j^n(t):=G_{\lambda_0}(t)+\int_{a}^{t}  P_{n}(K_{\lambda_0} (t,s,X_{j-1}^n(s))ds, \qquad (t \in [a,b]),
$$
where $P_n$ is described in Proposition \ref{pr:schauder}  using  as $\{f_n\}$ the usual basis in $C([0,1]^2)$ constructed  over the diadic nodes (see for instance \cite{Berenguer2011}). Fixing $\varepsilon >0$ we consider $m\in \mathbb{N}$ in such a way that
$$
H(X_m^n,X_{m-1}^n)< \varepsilon,
$$
and take $\tilde{X}(t):=X_m^n(t)$ as the target element. We can now address the inverse problem. Given the target element
 $\tilde{X}:[a,b]\rightarrow \I$, and fixed $r\in \N$, we consider 
\[
Y_{\lambda,r}(t)=G_\lambda(t)+\int_{a}^{t}  P_{r}(K_\lambda (t,s,\tilde{X}(s))ds, \qquad (t \in [a,b]).
\]
 We then compute  $\lambda_{r}^* \in \Lambda$ in such a way that
\begin{equation}\label{eq:min}
H(\tilde{X},Y_{\lambda_{r}^*,r})=\min_{\lambda\in\Lambda} H(\tilde{X},Y_{\lambda,r}), 
\end{equation}
and we analyse the differences between $\lambda_0$ and $\lambda_r^*$.

\begin{example}\label{ex:1}   Let 
\begin{equation} \label{eq:volterraej}
X(t)=G(t)+\int_{0}^{t} K(t,s,X(s))ds, \qquad (t \in [0,1])
\end{equation}
be the Volterra interval equation where $G(t)$ is obtained in such away that the solution of the equation is $X(t)=[\cos(t)-\frac{t}{2},\cos(t)+\frac{t}{2}]$ and $K(t,s,u)=(\sqrt{2} t-1 s)u$.  We now consider the family of Volterra interval integral equations
\[
X(t)=G_{\lambda}(t)+\int_{0}^{t} K_\lambda(t,s,X(s))ds, \qquad (t \in [0,1])
\]
with $G_\lambda(t)=G(t)$ and  $K_\lambda (t,s,u)=(\alpha t + \beta  s) u$, where $\lambda=(\alpha ,\beta ) $ with $ \alpha\in[1,3]$ and $\beta\in[-\frac{3}{2},-\frac{1}{2}].$ For the value $\lambda_0=(\alpha_0,\beta_0)=(\sqrt{2},-1)$, we compute $\tilde{X}=X_m^n$   and with this target element we  calculate $Y_{\lambda,r}$. Then we solve the minimization problem (\ref{eq:min}) and we note $\lambda_{r}^*=(\alpha_{r}^*, \beta_{r}^*)$ as the solution.  In Table 1, we show the obtained values for different $m$, $n$ and $r$.

\begin{table}[h]
\begin{center}

\small{\textsc{Table 1. Numerical results for Example \ref{ex:1}.}}

\begin{tabular}{|c|c|c|c|c|c|}
 \noalign{\hrule height 2pt} 
$m$  & $n$ &   $r$ & $\alpha_r^*$   &   $\beta_r^*$   & $H(\tilde{X},Y_{\lambda_{r}^*,r})$ \\  \noalign{\hrule height 2pt}
$3$  &  $3^2$  &  $3^2$  & $ 1.4052404924510171  $   & $ -0.9908622806761388 $  & $1.48128\times 10^{-9}$ \\  \hline
$3$  & $9^2 $ &  $9^2$  & $1.4127066294461754  $   & $-0.9982113290080197  $  & $3.14213\times 10^{-9}$ \\  \hline
$3$  & $17^2$  &  $17^2$  & $ 1.4064981933829088  $   & $-0.9904868353840546  $  &  $1.88601\times 10^{-9}$ \\  \hline
$7$  & $3^2 $ &  $3^2$  & $1.4142065762008067   $   & $-0.9999923158306506  $  & $2.58755\times 10^{-9}$  \\  \hline
$7$  & $9^2 $ &  $9^2$  & $1.4142138642529456  $   & $ -1.0000006138203197$  & $2.03712\times 10^{-9}$ \\  \hline
$7$  & $17^2 $ &  $17^2$  & $1.414213609421653   $   & $ -1.0000001074746905 $  & $2.8527\times 10^{-9}$ \\   \noalign{\hrule height 2pt}
\end{tabular}
\end{center}
\end{table}

\end{example}

\begin{example}\label{ex:2} 

We now consider the family of Volterra interval integral equations
\begin{equation} \label{eq:volterraej3}
X(t)=[2t+\frac{1}{8},2t+\frac{3}{8}]+\int_{0}^{t} (\alpha \cos(t)+\beta \cos(s)) \textrm{arctan}(X(s)) ds, \qquad (t \in [0,1])
\end{equation}
with $ \textrm{arctan}(X(s))=[\textrm{arctan}(\underline{X}(s)), \textrm{arctan}(\overline{X}(s))]$, $ \alpha\in[1.5,2.5]$ and $\beta\in[0.5,1.5]$. For $\alpha_0=2$ and $\beta_0=1$, we compute $\tilde{X}=X_m^n$   and  we  calculate $Y_{(\alpha,\beta),r}$. Then we solve the minimization problem (\ref{eq:min}) and we note $(\alpha_{r}^*, \beta_{r}^*)$ as the solution.

\begin{table}[h]
\begin{center}

\small{\textsc{Table 2. Numerical results for Example \ref{ex:2}.}}

\begin{tabular}{|c|c|c|c|c|c|}
 \noalign{\hrule height 2pt} 
$m$  & $n$ &   $r$ & $\alpha_r^*$   &   $\beta_r^*$   & $H(\tilde{X},Y_{\lambda_{r}^*,r})$ \\  \noalign{\hrule height 2pt}

$7$  & $9^2 $ &  $9^2$  & $1.9960170992122808  $   & $ 1.002917299482817$  & $3.33067\times 10^{-16}$ \\  \hline
$7$  & $17^2 $ &  $17^2$  & $1.9978735611639642  $   & $ 1.001405265860099 $  & $6.98024\times 10^{-9}$ \\   \noalign{\hrule height 2pt}
\end{tabular}
\end{center}
\end{table}

\end{example}

\section*{Conclusions}

In this paper we have established a perturbed Collage theorem. To deal with Volterra interval-valued integral equations, the use of adequate Schauder bases in certain Banach spaces, naturally associated with the inverse problem, allows us to design an algorithm based on the perturbed Collage Theorem.

In our future work, we will explore new contexts and applications related to this manuscript, which includes the vagueness of data needs of fuzzy modelling, as well as the necessary tools in order to extend the content and results of the present paper.

\section*{Acknowledgements}

\noindent Research partially supported by  MTM2017-
89577-P (MINECO, Spain), by project MTM2016-80676-P (AEI/Feder, UE), by Junta de Andaluc\'{i}a
Grant FQM359 and by E.T.S. Ingenier\'{i}a de Edificaci\'{o}n of
the University of Granada (Spain).


\begin{thebibliography}{99}

\bibitem{alefeld83} 
 \textsc{G. Alefeld, J. Herzberger}, 
{\sl Introduction to Interval Computations},
 Academic Press, New York, 1983.

\bibitem{alefeld-mayer10}
\textsc{ G. Alefeld, G. Mayer},
{\sl Interval analysis: theory and applications},
\textrm{ J. Comput. Appl. Math.} 
\textbf{121}  (2000) 421-464.

\bibitem{alm-ans-kha14}  
\textsc{S. Almezel, Q. H. Ansari. M. A. Khamsi} (eds.), 
{\sl Topics in fixed point theory}, 
Springer, Cham, 2014. 

\bibitem{Aubin-Cellina84}
 \textsc{J. P. Aubin, A. Cellina}, 
{\sl Differential Inclusions},
 Springer, New York, 1984.


\bibitem{Ervin}
\textsc{M. F.  Barnsley, V.  Ervin, D.  Hardin, J.  Lancaster },
{\sl  Solution of an inverse problem for fractals and other sets},
\textrm{ Proc. Nat. Acad. Sci. U.S.A.} \textbf{ 83} no. 7  (1986)  1975--1977. 



\bibitem{Bede-Gal05}  
\textsc{B. Bede, S. G. Gal}, 
{\sl  Generalizations of the differentiability of fuzzy number valued functions with applications to fuzzy differential equation},
 {\rm Fuzzy Sets Syst.} \textbf{151} (2005) 581--599.

\bibitem{Berenguer2017}
\textsc{M. I. Berenguer, D.  G\'amez }
{\sl  A computational method for solving a class of two dimensional Volterra integral equations},
{\rm J. Comput. Appl. Math. } \textbf{ 318} (2017), 403--410.

\bibitem{Berenguer2016}    
 \textsc{M. I. Berenguer, H. Kunze, D. La Torre, M. Ruiz Gal\'an}, 
{\sl  Galerkin method for constrained variational equations and a collage-based approach to related inverse problems.}, 
 {\rm J. Comput. Appl. Math.}  \textbf{292} (2016) 67--75.




\bibitem{Berenguer2011}  
\textsc{M. I. Berenguer, D.  G\'amez, A. I. Garralda-Guillem, M. Ruiz Gal\'an, M. C. Serrano P\'erez}, 
{\sl  Biorthogonal systems for solving Volterra integral equation systems of the second kind}, 
{\rm J. Comput. Appl. Math. } \textbf{23 }  no. 7 (2011) 1875--1883. 


\bibitem{Berenguer2010}  
 \textsc{M. I. Berenguer, D. G\'amez, A. I. Garralda-Guillem, M. C. Serrano P\'erez},
 {\sl Nonlinear Volterra integral equation of the second kind and biorthogonal systems}, 
  {\rm Abstr. Appl. Anal.} \textbf{2010}, Art. ID 135216, 11 pp.

\bibitem{Berenguer2009}    
 \textsc{M. I. Berenguer, D. G\'amez, A. I. Garralda-Guillem, M. Ruiz Gal\'an, M. C.  Serrano P\'erez}, 
{\sl  Analytical techniques for a numerical solution of the linear Volterra integral equation of the second kind}, 
 {\rm Abstr. Appl. Anal.}  \textbf{2009}, Art. ID 149367, 12 pp.

\bibitem{gelbaum}
  {\sc B. R. Gelbaum, J. Gil de Lamadrid},
{\sl  Bases of tensor products of Banach spaces},
  {\rm Pacific. J. Math.} {\bf 11} (1961) 1281--1286.



\bibitem{Kunze} 
\textsc{H. E. Kunze, E. R. Vrscay},
{\sl  Solving inverse problems for ordinary differential equations using the Picard contraction mapping}
{\rm Inverse Problems}, \textbf{15} (1999), 745--770.

\bibitem{KunzeTorre} 
\textsc{H. E. Kunze, D. La Torre, E. R. Vrscay},
{\sl  A generalized collage method upon the Lax-Milgran funtional for solving boundary value inverse problems}
{\rm Nonlinear Anal.}, \textbf{71} (2009), 1337--1343.



\bibitem{markov79}  
\textsc{S. Markov}, 
{\sl  Calculus for interval functions of real variable},
 {\rm Computing}, \textbf{ 22}  (1979) 325--337.

\bibitem{moore66}  
\textsc{R. E. Moore}, 
{\sl  Interval Analysis Prentice-Hall},
 Englewood Cliffs, NJ, 1966.

\bibitem{moore79}  
\textsc{R. E. Moore}, 
{\sl  Method and Applications of Interval Analysis},
 SIAM, Philadelphia, 1979.
 
 \bibitem{flores-medar02} \textsc{ H. Rom\'an-Flores, M. Rojas-Medar}
{\sl   Embedding of level-continuous fuzzy sets on Banach spaces}
\textrm{Information Sciences } \textbf{ 144 } (2002) 227--247.


\bibitem{semadeni1982}
  {\sc Z.  Semadeni},
{\sl  Schauder Bases in Banach Spaces of Continuous Functions},
  Springer--Verlag, Berl\'in, 1982.

\bibitem{shwartz16}
 \textsc{M. Shwartz, J. F. Burgess, J.  Zhu},  
{\sl  A DEA based composite measure of quality and its associated data
uncertainty interval for health care provider profiling and pay-for-performance},
{\rm Eur. J. Oper. Res.} \textbf{253} (2016)  489--502.

\bibitem{Stefanini2010} 
\textsc{L. Stefanini},
{\sl  A generalization of Hukuhara difference and division for interval and fuzzy
arithmetic},
 {\rm  Fuzzy Sets and Systems} \textbf{161 } (2010) 1564--1584.

\bibitem{Stefanini2009} 
\textsc{L. Stefanini, B. Bede},
{\sl  Generalized Hukuhara differentiability of interval-valued functions and interval differential equations},
 {\rm Nonlinear Anal.} \textbf{ 71 } (2009) 1311--1328.

\bibitem{stefanini-arana19}
 \textsc{L. Stefanini, M.  Arana-Jim\'enez}, 
{\sl  Karush-Kuhn-Tucker conditions for interval and fuzzy optimization in several variables under total and directional generalized differentiability},
{\rm  Fuzzy Sets and Systems} \textbf{362} (2019) 1--34 .   

\bibitem{WuGong00}
 \textsc{C. Wu, Z.  Gong}, 
{\sl  On Henstock integrals of inteval-valued functions and fuzzy-valued functions}, 
{\rm   Fuzzy Sets and Systems } \textbf{115} (2000) 377--391.





\end{thebibliography}
\end{document}